\documentclass[12pt,a4paper]{article}

\usepackage{latexsym}
\usepackage{amsfonts}
\usepackage{amsmath}

\newcommand{\diver}{{{\rm{div}\,}}}
\newcommand{\TM}{{T\hspace{-.1mm}M}}
\newcommand{\TO}{T\hspace{-.1mm}\Omega}
\newcommand{\OM}{{\Omega}}
\newcommand{\POM}{{\partial\Omega}}
\newcommand{\X}{{\cal X}}
\newcommand{\End}{{{\rm{End}\,}}}
\newcommand{\grad}{{{\rm grad\,}}}

\newcommand{\trace}{{{\rm Trace}\,}}
\newcommand{\hessian}{{{\rm Hess}\,}}
\newcommand{\noin}{\noindent}

\newtheorem{theorem}{\bf Theorem}[section]

\newtheorem{lemma}[theorem]{\bf Lemma}
\newtheorem{corollary}[theorem]{\bf Corollary}

\newtheorem{remark}[theorem]{\bf Remark}
\newtheorem{definition}[theorem]{\bf Definition}

\input amssym.def

\begin{document}

\title{Fundamental tone estimates for elliptic operators in divergence form and geometric applications }
\author{G. Pacelli Bessa
\and L. P. Jorge \\\and B. Pessoa Lima \and J. F\'{a}bio Montenegro}
\date{\today}
\maketitle
\begin{abstract} We establish a method for giving lower bounds for
the fundamental tone of  elliptic operators in divergence form in
terms of the divergence of vector fields. We then apply this method
to the $L_{r}$ operator associated to immersed hypersurfaces with
locally bounded $(r+1)$-th mean curvature $H_{r+1}$ of the space
forms $\mathbb{N}^{n+1}(c)$ of curvature $c$. As a corollary we give
lower bounds for the extrinsic radius  of closed hypersurfaces of
$\mathbb{N}^{n+1}(c)$ with $H_{r+1}>0$ in terms of the $r$-th and
$(r+1)$-th mean curvatures. Finally we observe that bounds for the
Laplace eigenvalues  essentially
  bound the eigenvalues of a
 self-adjoint elliptic differential operator in divergence form.  This allows us to
 show that Cheeger's constant gives a lower bounds for the first
 nonzero $L_{r}$-eigenvalue of a closed hypersurface of
 $\mathbb{N}^{n+1}(c)$.

\vspace{.2cm}
 \noindent
{\bf Mathematics Subject Classification } (2000): 58C40, 53C42

\vspace{.2cm}
 \noindent {\bf Key words:} Fundamental tone, closed and Dirichlet eigenvalue problems,  $L_{r}$ operator, $r$-th mean
 curvature, extrinsic radius, Cheeger's constant.

\end{abstract}
\section{Introduction}Let $\Omega$ be a domain in a smooth Riemannian
manifold $M$
  and let $\Phi:\OM \rightarrow \End(\TO) $ be a smooth symmetric and positive definite
  section of the bundle of all endomorphisms of $\TO$.  Each such
  section
   $\Phi$  is associated to a  second order self-adjoint  elliptic operator
$L_{\Phi}(f)=\diver (\Phi \;\grad \,f)$, $f \in C^{2}(\OM)$ so that
when $\Phi$ is the identity section then  $L_{\Phi}=\triangle$,  the
Laplace operator. The $L_{\Phi}$-fundamental tone of $\OM$ is
defined by \begin{equation}\label{FundTone-eq}
\lambda^{L_{\Phi}}(\OM) = \inf\left\{
 \displaystyle\frac{\int_{\OM}\vert \Phi^{1/2}\grad
 f\vert^{2}}{\int_{\OM}f^{2}}\,;\,f\in
 C_{0}^{2}(\OM)\setminus\{0\}\right\}.
\end{equation} If $\OM$ is bounded with smooth boundary $\POM\neq \emptyset$, the $L_{\Phi}$-fundamental
tone of $\OM$ coincides with the first eigenvalue
$\lambda_{1}^{L_{\Phi}}(\OM)$ of the Dirichlet eigenvalue problem
$L_{\Phi}\,u + \lambda\,u=0$ on $ \OM $, with $u\vert\partial
\OM=0$, $u\in C^{2}(\OM)
 \cap C^{0}(\overline{\OM})\setminus\{0\}$. If $\OM$ is bounded with
$\POM=\emptyset$ then $ \lambda^{L_{\Phi}}(\OM)=0$. A basic problem
  is what  lower and upper bounds for the fundamental tone of
a given domain $\OM$ in a smooth Riemannian manifold can be obtained
in terms of  Riemannian invariants of $\OM$.  In the first part of
this paper we show that the  method for giving lower bounds for the
$\triangle$-fundamental tone established in \cite{bessa-montenegro2}
can be extended for self-adjoint elliptic operators $L_{\Phi}$. The
lower bounds for the $L_{\Phi}$-fundamental tone of a domain $\OM$
are given in terms of the divergence of vector fields. By carefully
choosing  a test vector field, we can obtain lower bounds for the
$L_{\Phi}$-fundamental tone in terms of geometric invariants. This
is done in Theorem (\ref{Barta-L-thm}). We consider an immersed
hypersurface $M$ into the $(n+1)$-dimensional simply connected
 space form $\mathbb{N}^{n+1}(c)$ of constant sectional curvature $c\in\{1,0,-1\}$
  with locally bounded $(r+1)$-th mean curvature and such that a certain
differential operator $L_{r}$, $r\in\{ 0,1,\ldots, n\}$ is elliptic,
see \cite{Reilly}. Then we give lower bounds for the
$L_{r}$-fundamental tone of domains $\OM\subset
\varphi^{-1}(B_{\mathbb{N}^{n+1}(c)}(p,R))$ in terms of the $r$-th
and $(r+1)$-th mean curvatures $H_{r}$ and $H_{r+1}$. This is done
in Theorem
 (\ref{Hr+1-loc-thm}). We then derive  from this estimates three geometric corollaries (\ref{Hr+1-loc-cor1},
 \ref{Hr+1-loc-cor2}, \ref{Hr+1-loc-cor3}) that should be viewed as
 an extension of Theorem 1 of \cite{jorge-xavier}. There are related
 results due to Fontenele-Silva \cite{fontenele-silva}.
  To finish the first part of the paper we consider immersed hypersurfaces $M$
into $\mathbb{N}^{n+1}(c)$ such that the operators $L_{r}$ and
$L_{s}$, $0\leq r,s\leq n$ are elliptic and we compare the $L_{r}$
and $L_{s}$ fundamental
 tones $\lambda^{L_{r}}(\OM)$, $\lambda^{L_{s}}(\OM)$
 of domains $\OM\subset M\subset \mathbb{N}^{n+1}(c)$. In the
  second part of the
 paper we make an observation (Theorem \ref{Observation-thm})
  on the first nonzero eigenvalues of closed hypersurfaces. It follows that in order to get
  bounds for the eigenvalues of a self-adjoint elliptic differential operator
  $L_{\Phi}$ we essentially need bounds for the Laplace operator
  eigenvalues.
  This allows us to
 use Cheeger's constant  to give lower bounds for the first
 nonzero $L_{r}$-eigenvalue of a closed hypersurface of
 $\mathbb{N}^{n+1}(c)$.

\section{$L_{\Phi}$-fundamental tone estimates} Our main
estimate is the following method for giving lower bounds for
$L_{\Phi}$-fundamental tone of arbitrary domains of Riemannian
manifolds. It  extends the version of Barta's theorem \cite{barta}
proved by Cheng-Yau in \cite{cheng-yau}. It is the same proof (with
proper modifications) of a generalization of Barta's theorem proved
in \cite{bessa-montenegro2}.

\begin{theorem}\label{Barta-L-thm}Let $\OM$ be a domain in a Riemannian manifold
$M$ and let $\Phi:\OM \rightarrow \End(\TO) $ be a smooth symmetric
and positive definite
  section of $\TO$. Then the $L_{\Phi}$-fundamental tone of $\OM$ has the
  following lower bound
\begin{equation}\label{Barta-L-thm-eq2}\lambda^{L_{\Phi}}(\OM)\geq \sup_{{\cal X}(\Omega)}
 \inf_{\Omega} \left[ \diver (\Phi X)- \vert \Phi^{1/2} X
 \vert^{2}\right].
\end{equation} If $\OM$ is bounded and with smooth   boundary $\POM\neq \emptyset$
  then we have equality in (\ref{Barta-L-thm-eq2}).
\begin{eqnarray}\label{Barta-L-thm-eq1}\lambda^{L_{\Phi}}(\OM)&  = &  \sup_{\X(\OM )}
 \inf_{\OM} \left[ \diver (\Phi X)- \vert \Phi^{1/2} X
 \vert^{2}\right].
\end{eqnarray}
Where  $\X(\OM)$ is the set of all smooth  vector fields on $\OM$.
\end{theorem}  \section{Geometric applications}
Let us  consider the
 linearized operator $L_{r}$ of the $(r+1)$-mean curvature $ H_{r+1}=
 S_{r+1}/(\begin{array}{c}n\\r+1\end{array})$
 arising from normal variations of a
 hypersurface $M$   immersed into the $(n+1)$-dimensional simply connected space
  form $\mathbb{N}^{n+1}(c)$ of constant
 sectional curvature $c\in\{1,0,-1\}$ where
  $S_{r+1}$ is the $(r+1)$-th elementary symmetric function
of the principal curvatures $k_{1},k_{2},\ldots ,k_{n}$. Recall that
the elementary symmetric function of the principal curvatures are
given by
\begin{equation}
S_{0}=1,\;\;\;S_{r}= \sum_{i_{1}<\cdots <i_{r}} k_{i_{1}} \cdots
k_{i_{r}},\;\;1\leq r\leq n. \label{Sr-eq1}
\end{equation} Letting $A=-(\overline{\nabla}\eta)$ be  the shape operator
 of $M$, where $\overline{\nabla}$ is
the Levi-Civita connection of $\mathbb{N}^{n+1}(c)$ and $\eta$ a
globally defined unit vector field normal to $M$, we can recursively
define smooth symmetric sections $P_{r}:M\to \End(\TM)$, for $r =
0,1,\ldots, n$, called the Newton operators, setting
$P_{0}=I\,\,\,{\rm and }\,\, P_{r} = S_{r}Id-AP_{r-1}$ so that
$P_{r}(x): T_{x}M\rightarrow T_{x}M$ is a self-adjoint linear
operator  with the same eigenvectors as the shape operator $A$.
 \noin The operator $L_{r}$ is the
second order self-adjoint differential operator
\begin{equation}L_{P_{r}}(f)=\diver (P_{r}\, \grad f)\label{Lr-eq1}\end{equation} associated to the
section $P_{r}$. However, the sections $P_{r}$ may be not positive
definite and then the operators $L_{r}$ may not be elliptic, see
 \cite{Reilly}. However, there are geometric hypothesis that imply the ellipticity of $L_{r}$,
  see \cite{caffarelli-nirenberg-spruck}, \cite{korevaar}, \cite{barbosa-colares}. Here we will not impose
   geometric conditions to guarantee ellipticity of the $L_{r}$, except in corollary (\ref{Hr+1-loc-cor2}).
   Instead we will ask
   the ellipticity on the set of hypothesis in the following way. It is known, see
 \cite{kato},
 that there is an open and dense
   subset $U\subset M$
    where the ordered eigenvalues $\{\mu_{1}^{r}(x)\leq \ldots\leq
 \mu_{n}^{r}(x)\}$ of $P_{r} (x)$ depend smoothly on $x\in U$ and continuously on  $x\in M$. In
 addition, the respective eigenvectors $\{e_{1}(x),\ldots, e_{n}(x)\}$
form a smooth orthonormal frame in a neighborhood of
 every point of $U$.
Set $ \nu (P_{r})=\sup_{x\in M}\{ \mu_{n}^{r}(x)\}$
 and $ \mu
(P_{r})=\inf_{x\in M}\{ \mu_{1}^{r}(x)\} $. Observe that if $\mu
(P_{r})> 0$ then $P_{r}$ is positive definite, thus  $L_{r}$ is
elliptic.

 \vspace{2mm} \noin We need the following definition of
locally bounded $(r+1)$-th mean curvature hypersurface in order to
state our next result.
\begin{definition} An oriented  immersed
hypersurface $\varphi : M \hookrightarrow N$ of a Riemannian
manifold
 $N$ is said to have  locally bounded $(r+1)$-th mean curvature $H_{r+1}$
if  for any $p\in N$ and $R>0$, the number $h_{r+1}(p,R)=\sup
\{\vert S_{r+1}(x)\vert=a(n,r+1)\cdot \vert H_{r+1}(x)\vert
\,;\,x\in \varphi (M)\cap B_{N}(p,R) \} $ is finite.  Here  $
B_{N}(p,R)\subset N$ is the geodesic ball of radius $R$ and center
$p\in N$.
\end{definition}
Our next result generalizes in some aspects the main application of
\cite{bessa-montenegro}. There the first and fourth authors give
lower bounds for $\triangle$-fundamental tone of domains in
submanifolds with locally bounded mean curvature in complete
Riemannian manifolds.

\begin{theorem}
Let $\varphi \hspace{-.5mm}:M \hookrightarrow \mathbb{N}^{\,n+1}(c)$
be an oriented hypersurface immersed with locally bounded $(r+1)$-th
mean curvature $H_{r+1}$ for some $r\leq n-1$ and with $\mu
(P_{r})>0$. Let $B_{\mathbb{N}^{\,n+1}(c)}(p,R) $ be the geodesic
ball centered at $p\in \mathbb{N}^{\,n+1}(c)$ with radius $R$
 and $\OM\subset \varphi^{-1}(\,
\overline{B_{\mathbb{N}^{\,n+1}(c)}(p,R)}\,)$ be a connected
component.
 \label{Hr+1-loc-thm}Then the
$L_{r}$-fundamental tone $\lambda^{L_{r}}(\Omega)$ of $\OM$ has the
following lower bounds.

 \begin{itemize}\item[i.] For $c=1$ and $0<R< \displaystyle
   \cot^{-1}\left[\displaystyle\frac{(r+1)\cdot
   h_{r+1}(p,R)}{(n-r)\cdot \inf_{\Omega}
    S_{r}}\right]$ we have that

\begin{equation}\label{loc-eq1}\lambda^{L_{r}}(\Omega)\geq
   \displaystyle 2\cdot \frac{1}{R}
 \left[(n-r)\cdot\cot[ R]\cdot\inf_{\Omega}S_{r}-(r+1)\cdot  h_{r+1}(p,
  R)\right].\end{equation}
 \item [ii.] For $c\leq 0$, $ h_{r+1}(p,R)\neq  0$ and $0< R <  \displaystyle\frac{(n-r)\cdot
 \inf_{\Omega}
 S_{r}}{(r+1)\cdot  h_{r+1}(p,R)}$ we have that
\begin{equation}\label{loc-eq2}\lambda^{L_{r}}(\Omega)\geq
   2\cdot \displaystyle\frac{1}{R^{\,2}}
 \left[(n-r)\cdot\inf_{\Omega}S_{r}-(r+1)\cdot R\cdot h_{r+1}(p,
  R)\right].\end{equation}
  \item[iii.] If $c\leq 0$, $ h_{r+1}(p,R)=
 0$ and $R>0$ we have that
 \begin{equation}\label{loc-eq3}\lambda^{L_{r}}(\Omega)\geq
   \displaystyle\frac{2(n-r)\inf_{\Omega}S_{r}}{R^{\,2}}
 \end{equation}
\end{itemize}
\end{theorem}

 \begin{definition}Let $\varphi :M\hookrightarrow N$ be an isometric
 immersion of a closed   Riemannian manifold into a complete
 Riemannian manifold $N$. For each $x\in N$, let
$r(x)=\sup_{y\in M}
 dist_{N}(x,\varphi(y))$. The extrinsic radius $R_{e}(M)$ of $M$
 is defined by $$ R_{e}(M)=\inf_{x\in N}r(x).$$ Moreover, there is a
 point  $x_{0}\in N$ called the barycenter of
 $\varphi(M)$ in $N$
 such that
 $R_{e}(M)=r(x_{0})$. \end{definition}
  \begin{corollary} \label{Hr+1-loc-cor1}Let
  $\varphi :M \hookrightarrow B_{\mathbb{N}^{\,n+1}(c)}(R)\subset \mathbb{N}^{n+1}(c)$  be
  a complete
oriented hypersurface with  bounded $(r+1)$-th mean curvature
$H_{r+1}$ for some $r\leq n-1$,  $R$  chosen as in Theorem
(\ref{Hr+1-loc-thm}). Suppose that $\mu (P_{r})>0$ so that the
$L_{r}$ operator is elliptic. Then  $M$ is not closed.
  \end{corollary}
 \begin{corollary}\label{Hr+1-loc-cor2}Let $\varphi :M \hookrightarrow \mathbb{N}^{\,n+1}(c)$\footnote{If $c=1$ suppose
  that $\mathbb{N}^{n+1}(c)$ is the open hemisphere of $\mathbb{S}^{n+1}_{+}$.}, $c\in
 \{1,0,-1\}$  be an
 oriented closed hypersurface with $H_{r+1}>0$.
  Then there is an explicit constant $ \Lambda_{r} = \Lambda_{r}(c,\inf_{M}S_{r},\, \sup_{M} S_{r+1})>0$
   such that the extrinsic radius $R_{e}(M)\geq \Lambda_{r}$.
 \begin{itemize}\item [i.]For $c=1$,  $\Lambda_{r}= \displaystyle
   \cot^{-1}\left[\displaystyle\frac{(r+1)\cdot
   \sup_{M} S_{r+1}}{(n-r)\cdot \inf_{\Omega}
    S_{r}}\right] $.
    \item[ii.]For $c\in\{0,-1\}$,  $\Lambda_{r}= \displaystyle\frac{(n-r)\cdot
 \inf_{\Omega}
 S_{r}}{(r+1)\cdot \sup_{\OM} S_{r+1}}$.

 \end{itemize}

 \end{corollary}

\begin{remark}The hypothesis $H_{r+1}$ implies that $H_{j}>0$
and $L_{j}$ are elliptic for $j=0,1,\ldots r$,  see
\cite{barbosa-colares}, \cite{caffarelli-nirenberg-spruck} or
\cite{korevaar}. Thus  in fact in fact have that $R_{e}\geq
\max\{\Lambda_{0},\cdots, \Lambda_{r}\}$.
\end{remark}
\begin{remark}Jorge and Xavier, (Theorem 1 of \cite{jorge-xavier}), proved the inequalities of
 Corollary (\ref{Hr+1-loc-cor2}) when $r=0$
for complete submanifolds with
 scalar curvature bounded from below contained in a compact ball of a complete Riemannian manifold.
  Moreover, for $c=-1$ their
 inequality is
  slightly better. These inequalities
  should be also compared with a related result
  proved by
    Fontenele-Silva   in \cite{fontenele-silva}.
\end{remark}
 \begin{corollary}\label{Hr+1-loc-cor3} Let
   $\varphi :M \hookrightarrow \mathbb{S}^{\,n+1}(1)$,  be an
 oriented closed hypersurface with $\mu_{1}^{r}(M)>0$ and
 $H_{r+1}=0$. Then the extrinsic radius $R_{e}(M)\geq \pi/2$.
 \end{corollary}
 \begin{remark} An interesting question is:   Is it true that any
 closed oriented hypersurface with  $\mu_{1}^{r}(M)>0$ and
 $H_{r+1}=0$ intersect every great circle? For $r=0$  it is true and it  was proved
 by T.
 Frankel \cite{frankel}. \end{remark} We now consider immersed hypersurfaces
  $\varphi:M\hookrightarrow \mathbb{N}^{n+1}(c)$ with $L_{r}$ and
  $L_{s}$ elliptic. We can compare  the $L_{r}$ and $L_{s}$
  fundamental tones of a domain $\OM\subset M$. In particular we can
  compare with its $L_{0}$-fundamental tone.
\begin{theorem}\label{thm1.5}Let $\varphi:M\hookrightarrow \mathbb{N}^{n+1}(c)$ be an
 oriented $n$-dimensional hypersurface
$M$ immersed into the  $(n+1)$-dimensional simply connected space
form of constant sectional curvature $c$ and $\mu(L_{r})>0$ and
$\mu(L_{s})>0$, $0\leq s,\, r\leq n-1$. Let $\OM \subset M$ be a
domain with compact closure and piecewise smooth non-empty boundary.
Then the $L_{r}$ and $L_{s}$ fundamental tones satisfies the
following inequalities
 \begin{equation}\label{eqPrPs}
 \lambda^{L_{r}}(\OM)\geq \displaystyle\frac{\mu(P_{r})}{\nu(P_{s})}\cdot\lambda^{L_{s}}(\OM)
 \end{equation}Where  $ \lambda^{L_{s}}(\OM)$ and  $ \lambda^{L_{r}}(\OM)$ are respectively the first
 $L_{s}$-eigenvalue and $L_{r}$-eigenvalue
 of $\OM$. From (\ref{eqPrPs})  we have in particular  that
 \begin{equation}\label{eqPrP0}
 \nu(r)\cdot\lambda^{\triangle}(\OM)\geq
 \lambda^{L_{r}}(\OM)\geq
 \mu(r)\cdot\lambda^{\triangle}(\OM)
 \end{equation}

 \end{theorem}

 \subsection{Closed eigenvalue problem}
Let $M$ be a closed hypersurface of a simply connected space form
$\mathbb{N}^{n+1}(c)$.  Similarly to the eigenvalue problem of
closed Riemannian manifolds, the interesting  problem is what bounds
can one obtain for the first nonzero $L_{r}$-eigenvalue
$\lambda_{1}^{L_{r}}(M)$ in terms of the geometries of $M$ and of
the ambient space.
 Upper bounds for the first nonzero $\triangle$-eigenvalue or even for the first nonzero
 $L_{r}$-eigenvalue, $r\geq 1$
 have been obtained  by many authors
 in contrast with  lower bounds that are rare. For instance, Reilly \cite{Reilly2} extending earlier result of
Bleecker and Weiner \cite{bleecker-weiner} obtained  upper bounds
for $\lambda_{1}^{\triangle}(M)$ of a closed submanifold $M$ of
$\mathbb{R}^{m}$ in terms of the total mean curvature of $M$.
Reilly's result applied to compact submanifolds of the sphere
$M\subset\mathbb{S}^{m+1}(1)$, this later viewed as a hypersurface
of the Euclidean space $\mathbb{S}^{m+1}(1)\subset \mathbb{R}^{m+2}$
obtains upper bounds for $\lambda_{1}^{\triangle}(M)$, see
\cite{alencar-rosenberg}. Heintze,\cite{heintze} extended Reilly's
result to  compact manifolds
 and  Hadamard manifolds $\overline{M}$. In particular for the hyperbolic space $\mathbb{H}^{n+1}$.
 The best upper bounds for the
 first nonzero
$\triangle$-eigenvalue of   closed hypersurfaces $M$ of
$\mathbb{H}^{n+1}$ in terms of the total mean curvature of $M$ was
obtained by El Soufi and Ilias \cite{soufi-ilias}. Regarding the
$L_{r}$ operators, Alencar, Do Carmo, and  Rosenberg
  \cite{alencar-rosenberg} obtained   sharp (extrinsic) upper bound
  the first nonzero
eigenvalue  $\lambda_{1}^{L_{r}}(M)$ of the linearized
 operator $L_{r}$ of  compact hypersurfaces $M$ of  $\mathbb{R}^{m+1}$ with  $S_{r+1}>0$.
  Upper bounds for $\lambda_{1}^{L_{r}}(M)$ of compact
 hypersurfaces of $\mathbb{S}^{n+1}$, $\mathbb{H}^{n+1}$ under the hypothesis that $L_{r}$ is elliptic were obtained
 by Alencar, Do Carmo, Marques in \cite{alencar-coda} and by Alias and Malacarne in \cite{alias-malacarne}
 see also the work of Veeravalli \cite{veeravalli}.
On the other hand,  lower
  bounds for $\lambda_{1}^{L_{r}}(M)$ of closed
 hypersurfaces $M\subset\mathbb{N}^{n+1}(c)$ are not so well studied
as the  upper bounds, except for $r=0$ in which case
$L_{0}=\triangle$. In this paper we make a simple observation
(Theorem \ref{Observation-thm})  that to obtain lower and upper
bounds for the
  $L_{\Phi}$-eigenvalues (Dirichlet or Closed eigenvalue problem)
  it is enough to obtain
  lower and upper bounds for
   the eigenvalues of $\Phi$ and for the
eigenvalues for the Laplacian in the respective problem. When
applied to the $L_{r}$ operators (supposing them elliptic) we obtain
lower bounds for closed hypersurfaces of the space forms via
Cheeger's lower bounds for the first $\triangle$-eigenvalue of
closed manifolds. Let   $\{\mu_{1}(x)\leq \ldots\leq
 \mu_{n}(x)\}$ be the ordered eigenvalues of $\Phi (x)$.  Setting $ \nu
(\Phi)=\sup_{x\in \OM}\{ \mu_{n}(x)\}$ and $ \mu (\Phi)=\inf_{x\in
\OM}\{ \mu_{1}(x)\} $ we have the following theorem.
\begin{theorem}\label{Observation-thm} Let $\lambda^{L_{\Phi}}(\OM)$ denote the
$L_{\Phi}$-fundamental tone of $\OM$ if $\OM$ is unbounded or
$\POM\neq \emptyset$ and the first nonzero $L_{\Phi}$-eigenvalue
$\lambda_{1}^{L_{\Phi}}(\OM)$ if  $\OM$  is a closed manifold. Then
 $\lambda^{L_{\Phi}}(\OM)$
  satisfies the following
 inequalities,
\begin{equation}\label{Obs-eq1}\nu(\Phi,\OM)\cdot\lambda^{\triangle}(\OM)\geq \lambda^{L_{\Phi}}(\OM)
\geq \mu(\Phi,\OM)\cdot \lambda^{\triangle}(\OM),\end{equation}
where $\lambda^{\triangle}(\OM) $ is the $\triangle$-fundamental
tone of $\OM$ or the first nonzero $\triangle$-eigenvalue of $\OM$.
\end{theorem}Let $M$ be a closed $n$-dimensional Riemannian manifold, Cheeger in
\cite{cheeger} defined  the following  constant given by
\begin{equation}\label{CheegerConst-eq}
h(M)=\inf_{S}\frac{vol_{n-1}(S)}{\min\{vol_{n}(\OM_{1}),\,vol_{n}(\OM_{2})\}},\end{equation}
  where $S\subset M$ ranges over all
connected closed hypersurfaces dividing $M$ in two connected
components, i.e.
 $M=\OM_{1}\cup\OM_{2}$, $\OM_{1}\cap\OM_{2}=\emptyset$ such that $S=\POM_{1}=\POM_{2}$ and he proved that
the first nonzero $\triangle$-eigenvalue
$\lambda_{1}^{\triangle}(M)\geq h(M)^{2}/4$.
\begin{corollary}Let $\varphi :M \hookrightarrow \mathbb{N}^{\,n+1}(c)$, $c\in
 \{1,0,-1\}$\footnote{If $c=1$ suppose
  that $\mathbb{N}^{n+1}(c)$ is the open hemisphere of $\mathbb{S}^{n+1}_{+}$.}  be an
 oriented closed hypersurface with $H_{r+1}>0$. Then the first nonzero $L_{r}$-eigenvalue of $M$
 has the following lower bound $$ \lambda_{1}^{L_{r}}(M)
\geq \mu(L_{r})\cdot \frac{h^{2}(M)}{4}
$$
\end{corollary}

\section{Proof of the Results}
\subsection{Proof of Theorem \ref{Barta-L-thm}}
 Let $\Omega$ be an arbitrary domain, $X$ be  a smooth vector field on $ \Omega $ and
 $f\in C^{\infty}_{0}(\Omega)$. The vector field $f^{\,2}\Phi X$ has compact support
 ${\rm supp}(f^{\,2}\Phi X)\subset {\rm supp}(f) \subset\Omega
 $. Let $ {\cal S}$ be a  regular domain containing the support of
 $f$.
 We have by the divergence theorem that
  \begin{eqnarray}0=\int_{{\cal S}}\diver (f^{2}\Phi X) & =   & \int_{\Omega}\diver (f^{2}\Phi X) \nonumber \\
   && \nonumber \\
   & = & \int_{\Omega}\left[\langle \grad f^{2}, \Phi X\rangle + f^{2}\diver (\Phi X)\right]\nonumber \\
   & & \nonumber \\
                            &\geq & -2 \, \int_{\Omega}\left[\vert f\vert \cdot \vert \Phi ^{1/2}\grad f\vert\cdot
                             \vert \Phi^{1/2} X \vert+
                             \diver (\Phi X)\cdot
                            f^{2}\right]
                             \label{ProofThm1.1a} \\
                             && \nonumber \\
                             & \geq &  \, \int_{\Omega}\left[- \vert \Phi ^{1/2}\grad f\vert^{2} - f^{2} \cdot
                                                          \vert \Phi^{1/2} X \vert^{2}+
                             \diver (\Phi X)\cdot
                            f^{2}\right].\nonumber
\end{eqnarray}

Therefore
\begin{eqnarray} \label{ProofThm1.1b}
  \int_{\Omega}\vert \Phi^{1/2} \grad f\vert^{2}&\geq & \int_{\Omega}\left[ \diver (\Phi X)- \vert \Phi^{1/2} X \vert^{2}
  \right] f^{2}\nonumber\\
  & \geq & \inf \left[ \diver (\Phi X)- \vert \Phi^{1/2} X \vert^{2}\right] \int_{\Omega}f^{2}
\end{eqnarray} By
the variational formulation (\ref{FundTone-eq}) of
$\lambda^{L_{r}}(\Omega)$ this inequality above implies that
 \begin{equation}\lambda^{L_{r}}(\Omega)\geq  \inf_{\Omega} \left[ \diver (\Phi X)- \vert \Phi^{1/2} X
 \vert^{2}\right].
  \label{ProofThm1.1c}\end{equation}When $\Omega$ is a bounded domain with smooth boundary
  $\POM \neq \emptyset$ then
  $\lambda^{L_{r}}(\Omega)= \lambda_{1}^{L_{r}}(\Omega)$. This proof above  shows
  that $\lambda^{L_{r}}_{1}(M) \geq  \inf_{M} \left[ \diver (\Phi X)- \vert \Phi^{1/2} X
 \vert^{2}\right]$. Let $v \in C^{2}(\OM)\cap C^{0}(\overline{\OM})$
  be a positive first $L_{r}$-eigenfunction\footnote{
 $v\in C^{2}(\OM)\cap H_{1}^{0}(\OM)$ if $\POM$ is not smooth.} of
 $\OM$ and if we set
   $X_{0}=-\grad \log (v)$ we have that

 \begin{equation}\label{ProofThm1.1d}\begin{array}{lcl}\diver (\Phi X_{0})- \vert \Phi^{1/2} X_{0}
 \vert^{2} &= & -\diver \,((1/v)\,\Phi \,\grad v )- (1/v^{2})\,\vert \Phi^{1/2}\,\grad v \vert^{2}\\
 && \\
 &= & (1/v^{2})\langle \grad v , \Phi \,\grad v \rangle - (1/v) \,\diver (\Phi\, \grad
 v)\\
 && \\
 && -(1/v^{2})\,\vert \Phi^{1/2}\,\grad v \vert^{2}\\
 && \\
 &=&- (1/v) \,\diver (\Phi\, \grad
 v)=-L_{r} (v)/v= \lambda_{1}^{L}(\OM).\\

 \end{array} \end{equation}This proves (\ref{Barta-L-thm-eq1}).

 \subsection{Proof of Theorem
\protect{\ref{Hr+1-loc-thm}} and   Corollaries \ref{Hr+1-loc-cor1},
\ref{Hr+1-loc-cor2}, \ref{Hr+1-loc-cor3}}We start this section
stating few lemmas necessary to construct the proof of Theorem
(\ref{Hr+1-loc-thm}). The first  lemma  was proved in
\cite{jorge-koutrofiotis} for the Laplace operator  and for the
$L_{r}$ operator in \cite{pessoa-lima1} and \cite{pessoa-lima2}. We
reproduce its proof to make the exposition complete.
 \begin{lemma}\label{lemma2.1}Let $\varphi :M\hookrightarrow
 \mathbb{N}^{n+1}(c)$  be a complete hypersurface immersed in $(n+1)$-dimensional simply
  connected space form $\mathbb{N}^{n+1}(c)$
 of constant sectional curvature $c$. Let
 $g:\mathbb{N}^{n+1}(c)\rightarrow\mathbb{R}$ be a smooth
 function and set $f=g\circ \varphi $. Identify $X \in T_{p}M$
  with $d\varphi (p) X \in T_{\varphi
 (p)}\,\varphi (M)$ then we
  have that
 \begin{equation}\label{lemma2.1a}\begin{array}{lll}
 L_{r}f (p)&=&\sum_{i=1}^{n}\mu_{i}^{r}\,\hessian  g(\varphi (p))\, (e_{i},e_{i})+\trace (A P_{r})\langle \grad g, \eta \rangle
 \end{array}
 \end{equation}
\end{lemma}
\noin {\em Proof:} Each  $P_{r}$  is also associated to a second
order self-adjoint differential operator defined by $\Box f=\trace
(P_{r}\, \hessian (f))$ see \cite{cheng-yau}, \cite{hartmann}. We
have that
\begin{equation}\Box f=\trace (P_{r}\, \hessian (f))=\diver (P_{r}\,\grad f)-{\rm trace}\left( \nabla
P_{r}\right)\grad f .\label{eqLrBox}\end{equation}Rosenberg
\cite{rosenberg} proved that when the ambient manifold   is the
simply connected space form $\mathbb{N}^{n+1}(c) $ then  ${\rm
Trace}\left( \nabla P_{r}\right)\,\grad \equiv 0 $, see also
\cite{Reilly}. Therefore $L_{r}f=\trace (P_{r}\, \hessian (f))$.
Using Gauss equation to compute $\hessian (f)$ we obtain
\begin{equation}\label{lemma2.1b} \hessian f(p)(X,Y)=\hessian g(\varphi (p))(X,Y) +
\langle \grad g, \alpha (X,Y)\rangle_{\varphi (p)},
\end{equation}where   $\langle \alpha
 (X,Y),\eta\rangle=\langle A(X),Y\rangle$.
Let $\{e_{i}\}$ be an orthonormal frame around $p$ that diagonalize
the section $P_{r}$ so that
 $P_{r}(x)(e_{i})=\mu_{i}^{r}(x)e_{i}$.
Thus
\begin{equation}\label{lemma2.1c}\begin{array}{lll}
   L_{r}f & = & \sum_{i=1}^{n}\langle P_{r}\,\hessian f
(e_{i}),e_{i}\rangle  \\
 && \\
   & = & \sum_{i=1}^{n}\langle \,\hessian f
(e_{i}),\mu_{i}^{r}e_{i}\rangle  \\
&& \\
   & = & \sum_{i=1}^{n}\mu_{i}^{r}
\,\hessian f (e_{i},e_{i})
\end{array}
\end{equation} Substituting (\ref{lemma2.1b}) into (\ref{lemma2.1c}) we have
that \begin{equation}\label{lemma2.1d}\begin{array}{lll}
    L_{r}f& = & \sum_{i=1}^{n}\mu_{i}^{r}
\,\hessian g \,(e_{i},e_{i})+ \langle \grad g,
\sum_{i=1}^{n}\mu_{i}^{r}\alpha (e_{i},e_{i})
\rangle  \\
&& \\
&=& \sum_{i=1}^{n}\mu_{i}^{r} \,\hessian g\, (e_{i},e_{i})+
\langle \grad g, \alpha (\sum_{i=1}^{n}P_{r}(e_{i}),e_{i})
\rangle    \\
&& \\
 & = & \sum_{i=1}^{n}\mu_{i}^{r} \,\hessian g \,(e_{i},e_{i})+
\trace (A P_{r})\langle \grad g, \eta \rangle
\end{array}
\end{equation}
Here $\hessian f (X)=\nabla_{X} \grad f$ and $\hessian f
(X,Y)=\langle \nabla_{X} \grad f,Y\rangle$. The next two lemmas we
are gong to present are well known and their proofs are easily found
in the literature thus we will omit them here.

\begin{lemma}[Hessian Comparison Theorem] Let $M$ be    a  complete Riemannian manifold  and
$x_{0},x_{1} \in M $.   Let
 $\gamma:[0,\,\rho (x_{1}) ]\rightarrow M$ be a minimizing geodesic joining $x_{0}$ and $x_{1}$ where $\rho (x)$
 is the  distance function $dist_{M}(x_{0}, x) $. Let $K$ be the sectional curvatures of $M$  and  $\upsilon(\rho )$,
  defined below.
\begin{equation} \upsilon(\rho )=\left\{ \begin{array}{lcll}
 & k_{1} \cdot\coth (k_{1} \cdot\rho (x)), & if  & \sup_{\gamma} K=-k_{1}^{2} \\
 &                       &     &   \\
 & \displaystyle\displaystyle\frac{1}{\rho (x)},   & if  &  \sup_{\gamma} K=0  \\
 &                       &     & \\
 & k_{1} \cdot\cot (k_{1}\cdot \rho (x)),  & if  &  \sup_{\gamma} K =k_{1}^{2}\; and \; \rho < \pi/2k_{1}.
\end{array}\right.\label{hessianA}
\end{equation}
\vspace{.1cm}

\noin Let $X=X^{\perp}+X^{T}\in T_{x}M$, $X^{T}=\langle
X,\gamma'\rangle \gamma'$ and $\langle X^{\perp},\,
\gamma'\rangle=0$. Then
\begin{equation}\label{hessianB}
Hess\,\rho(x)(X,X)= Hess\,\rho(x)(X^{\perp},X^{\perp})\geq
\upsilon(\rho(x))\cdot\Vert X^{\perp}\Vert^{2}
\end{equation}
\end{lemma} See \cite{schoen-yau} for a proof.
\begin{lemma}\label{Properties-lem}Let $p\in M$ and  $1\leq r\leq n-1$, let
    $\{e_{i}\}$ be an orthonormal basis of $T_{p}M$ such that
   $P_{r}(e_{i})=\mu_{i}^{r} e_{i}$ and
   $A(e_{i})=k_{i}e_{i}$. Then
   \begin{itemize}
   \item[i.] ${\rm trace}(P_{r})= \sum_{i=1}^{n}\mu_{i}^{r}=(n-r)S_{r}$
   \item[ii.]$ {\rm trace}(A P_{r})=\sum_{i=1}^{n}k_{i}\mu_{i}^{r}=(r+1)S_{r+1}$
   \end{itemize}In particular, if the Newton operator  $P_{r}$ is positive
   definite then
   $S_{r}>0$.
   \end{lemma}

\noin To prove Theorem (\ref{Hr+1-loc-thm})  set
 $g:B(p, R)\subset \mathbb{N}^{n+1}(c)\rightarrow \mathbb{R}$
  given by $g=R^{\,2}-\rho^{2} $, where $\rho $
  is the distance
function  ($\rho(x)={\rm dist}(x,p)$) of $\mathbb{N}^{n+1}(c)$.
Setting $f=g\circ\varphi $ we obtain by  (\ref{lemma2.1a}) that
\begin{equation}\label{hessianC}L_{r}f=\sum_{i=1}^{n}\mu_{i}^{r} \cdot\hessian g
\,(e_{i},e_{i})+(r+1)\cdot S_{r+1}\cdot\langle \grad g, \eta
\rangle,
\end{equation}since  $\trace (A P_{r})=(r+1)\cdot S_{r+1}$. Letting  $X =-\grad \log
f$ we have by Theorem (\ref{Barta-L-thm}) that
\begin{equation}\label{hessianD}
\lambda^{L_{r}}(\Omega) \geq
\inf_{\Omega}(-L_{r}f/f)=\inf_{\Omega}\left\{-\displaystyle\frac{1}{g}\left[
\sum_{i=1}^{n}\mu_{i}^{r} \cdot\hessian g\, (e_{i},e_{i})+(r+1)\cdot
S_{r+1}\cdot\langle \grad g, \eta \rangle\right]\right\}.
\end{equation}
  Computing the Hessian of $g$
we have that
\begin{equation}\label{hessianE}\begin{array}{lll}\hessian g \,(e_{i},e_{i})&= &\langle \nabla_{e_{i}}\grad
g, e_{i}\rangle = -2\langle \nabla_{e_{i}}\rho \, \grad
\rho, e_{i}\rangle\\
&& \\
&=& -2\langle \grad \rho, e_{i}\rangle^{2} - 2\rho \, \langle
\nabla_{e_{i}} \grad
\rho, e_{i}\rangle\\
&&\\
&=&-2\langle \grad \rho, e_{i}\rangle^{2}- 2\rho \,\hessian \rho
(e_{i},e_{i}).
\end{array}
\end{equation}Therefore we have that
\begin{equation}\label{hessianF}-\frac{L_{r}f}{f}=\frac{2}{R^{\,2}-\rho^{2}}
\left[\sum_{i=1}^{n}\mu_{i}^{r} \,[\langle \grad \rho,
e_{i}\rangle^{2}+ \rho \,\hessian \rho (e_{i},e_{i})]+(r+1)\cdot
S_{r+1}\cdot\rho\cdot \langle \grad \rho, \eta \rangle\right]
\end{equation}Setting $e_{i}^{T}=\langle \grad \rho,
e_{i}\rangle \grad \rho $ and $e_{i}^{\perp}=e_{i}-e_{i}^{T}$, by
the Hessian Comparison Theorem  we have that
\begin{equation}\label{hessianG}\sum_{i=1}^{n}\mu_{i}^{r} [\langle \grad \rho,
e_{i}\rangle^{2}+ \rho \hessian \rho (e_{i},e_{i})] \geq
\sum_{i=1}^{n}\mu_{i}^{r}\left[\Vert e_{i}^{T}\Vert^{2}+ \rho\cdot
\upsilon (\rho )
  \Vert e_{i}^{\perp}\Vert^{2}\right]
\end{equation} and \begin{equation}\label{hessianH}
(r+1)\cdot S_{r+1}\cdot\rho\cdot \langle \grad \rho, \eta \rangle
\leq (r+1)\,R\,\cdot h_{r+1}(p, R)
\end{equation}From (\ref{hessianG}) and (\ref{hessianH}) wee have that
\begin{eqnarray}\label{hessianI}\lambda^{1}(\Omega)&  \geq &
 \inf_{\Omega}(-L_{r}f/f)\nonumber \\
  & \geq & 2\cdot \inf_{\Omega}\left\{\frac{1}{R^{\,2}-\rho^{2}}
 \left[\sum_{i=1}^{n}\mu_{i}^{r}\left[\Vert e_{i}^{T}\Vert^{2}+
\rho\cdot \upsilon (\rho )
  \Vert e_{i}^{\perp}\Vert^{2}\right]-(r+1)\cdot R\,\cdot h_{r+1}(p,
  R)\right]\right\}
\end{eqnarray} If $c\leq 0$ then $ \rho\cdot \upsilon (\rho )\geq 1$ thus
from (\ref{hessianI})  we have that
\begin{eqnarray}\lambda^{1}(\Omega) & \geq & 2\cdot \frac{1}{R^{\,2}}
 \left[\inf_{\Omega}\left\{\sum_{i=1}^{n}\mu_{i}^{r}\left[\Vert e_{i}^{T}\Vert^{2}+
  \Vert e_{i}^{\perp}\Vert^{2}\right]\right\}-(r+1)\cdot R \cdot h_{r+1}(p,
  R)\right] \nonumber \\
  & = & 2\cdot \frac{1}{R^{\,2}}
 \left[\inf_{\Omega}\sum_{i=1}^{n}\mu_{i}^{r}-(r+1)\cdot R \cdot h_{r+1}(p,
  R)\right]\label{hessianJ} \\
  & =& 2\cdot \frac{1}{R^{\,2}}
 \left[(n-r)\inf_{\Omega}S_{r}-(r+1)\cdot R \cdot h_{r+1}(p,
  R)\right].\nonumber
\end{eqnarray}
If $c>0$ then $ \rho\cdot \upsilon (\rho
)=\rho\cdot\sqrt{c}\cdot\cot[\sqrt{c}\,\rho]\leq 1$  thus from
(\ref{hessianI}) we have that
\begin{eqnarray}\lambda^{1}(\Omega) & \geq & 2\cdot \frac{1}{R^{\,2}}
 \left[\inf_{\Omega}\left\{\sum_{i=1}^{n}\mu_{i}^{r}\left[\Vert e_{i}^{T}\Vert^{2}+
  \Vert e_{i}^{\perp}\Vert^{2}\right]\,\rho\cdot\sqrt{c}\cdot\cot[\sqrt{c}\,\rho]\right\}-(r+1)\cdot R\cdot h_{r+1}(p,
  R)\right]\nonumber \\
  & = & 2\cdot \frac{1}{R^{\,2}}
 \left[\inf_{\Omega}\left\{\sum_{i=1}^{n}\mu_{i}^{r}\,\rho\sqrt{c}\cot[\sqrt{c}\rho]\right\}-(r+1)\cdot R\cdot h_{r+1}(p,
  R)\right]\label{hessianL}\\
  & =& 2\cdot \frac{1}{R^{\,2}}
 \left[(n-r)\cdot R\cdot\sqrt{c}\cdot\cot[\sqrt{c}R]\cdot\inf_{\Omega}S_{r}-(r+1)\cdot R \cdot h_{r+1}(p,
  R)\right].\nonumber
\end{eqnarray}

\noin To prove the Corollaries (\ref{Hr+1-loc-cor1}) and
(\ref{Hr+1-loc-cor2}), observe that the hypotheses $\mu(P_{r})(M)>0$
(in Corollary \ref{Hr+1-loc-cor1}) and $H_{r+1}>0$  (in Corollary
\ref{Hr+1-loc-cor2}) imply that the $L_{r}$ is elliptic. If the
immersion is bounded (contained in a ball of radius $R$, for those
choices of $R$) and  $M$ is closed we would have by one hand that
the $L_{r}$-fundamental tone would be zero  and by Theorem
(\ref{Hr+1-loc-thm}) that it would be positive. Then $M$ can not be
closed if the immersion is bounded. On the other hand if $M$ is
closed a ball of radius $R$ centered at the barycenter of $M$ could
not contain $M$ because the fundamental tone estimates for any
connected component $\OM \subset \varphi^{-1}(\varphi (M)\cap
B_{\mathbb{N}^{n+1}(c)}(p,R)$  is positive. Showing that $M\neq
\OM$.
  The corollary
(\ref{Hr+1-loc-cor3}) follows from item i. of Theorem
(\ref{Hr+1-loc-thm}) placing $S_{r+1}=0$.

\subsection{Proof of  Theorem \ref{thm1.5}}Let $\varphi :W \hookrightarrow
\mathbb{N}^{n+1}(c)$ be an isometric immersion of an oriented
$n$-dimensional Riemannian manifold $W$ into a $(n+1)$-dimensional
simply connected space form of sectional curvature $c$. Let $M
\subset W$ be a  domain with compact closure and piecewise smooth
nonempty boundary
 and suppose that the Newton
operators  $P_{r}$ and $P_{s}$, $0\leq s,\,r\leq n-1$ are positive
definite when restricted to $M$. Let $\mu (r)=\mu (P_{r},M)$, $\mu
(s)=\mu (P_{s},M)$  and $\nu (r)=\nu (P_{r},M )$,   $\nu (s)=\nu
(P_{s},M )$.
 Given a vector field $X$ on $M $ we can find a vector field $Y$ on $M$  such that
  $P_{r}X=\kappa \cdot P_{s}Y$, $\kappa $ constant. Now
 \begin{eqnarray}\label{eq22}\diver (P_{r}X )-\vert
 P_{r}^{1/2}X\vert^{2} & = & \kappa \cdot \diver (P_{s}Y)-\langle
 P_{r}X, X\rangle\nonumber \\
  & = & \kappa \cdot \diver (P_{s}Y)-\kappa^{2}\langle P_{s}Y, P_{r}^{-1}P_{s}Y\rangle\\
  & = & \kappa \cdot\left[\diver
  (P_{s}Y)-\vert P_{s}^{1/2}Y\vert^{2}+\vert P_{s}^{1/2}Y\vert^{2}-
  \kappa \cdot \vert P_{r}^{-1/2}P_{s}Y\vert^{2}\right ]\nonumber
  \end{eqnarray}Consider $\{e_{i}\}$ be an orthonormal basis such
  that $P_{r}e_{i}=\mu_{i}^{r}e_{i}$ and
  $P_{s}e_{i}=\mu_{i}^{s}e_{i}$. Letting $Y=\sum_{i=1}^{n}y_{i}e_{i}$ then
  \begin{eqnarray}\vert P_{s}^{1/2}Y\vert^{2}-
  \kappa \cdot \vert P_{r}^{-1/2}P_{s}Y\vert^{2} &
  =&\sum_{i=1}^{n}\mu_{i}^{s}\,y_{i}^{2}-\kappa
  \sum_{i=1}^{n}\frac{(\mu_{i}^{s})^{2}}{\mu_{i}^{r}}\,y_{i}^{2}\nonumber
  \\ &= &\sum_{i=1}^{n}\mu_{i}^{s}\,y_{i}^{2}\left[ 1- \kappa
  \cdot \frac{\mu_{i}^{s}}{\mu_{i}^{r}}\right]\label{eq23}\\
  & \geq & 0, \,\,\, if \,\,\,  \kappa \leq   \displaystyle\frac{\mu (r)}{\nu
   (s)}\nonumber
  \end{eqnarray} Combining (\ref{eq22}) with (\ref{eq23}) and by
  Theorem (\ref{Barta-L-thm}) we have that
\begin{equation} \lambda^{L_{r}}(M)=\sup_{X}
\inf_{M}\diver (P_{r}X )-\vert P_{r}^{1/2}X\vert^{2}\geq \kappa
\cdot\sup_{Y}\inf_{M}\diver (P_{s}Y) -\vert
P_{s}^{1/2}Y\vert^{2}=\kappa \cdot \lambda^{L_{s}}(M),\label{eq24}
\end{equation} for every $0<\kappa \leq  \displaystyle\frac{\mu (r)}{\nu
   (s)}$. This proves (\ref{eqPrPs}).
 \subsection{Proof of  Theorem \ref{Observation-thm}}Recall that for any smooth symmetric section
 $\Phi:\OM\to \End (\TO)$ there is an open and dense subset $U\subset
\OM$ where the ordered eigenvalues $\{\mu_{1}(x)\leq \ldots\leq
 \mu_{n}(x)\}$ of $\Phi (x)$ depend smoothly on $x\in U$ and continuously in all $\OM$. In
 addition, the respective eigenvectors $\{e_{1}(x),\ldots,
 e_{n}(x)\}$ form a smooth orthonormal frame in a neighborhood of
 every point of $U$, see \cite{kato}.
  Let $f\in C_{0}^{2}(\OM)\setminus \{0\}$ ($f \in C^{2}(\OM)$ with $\smallint_{\OM}f=0$) be
an admissible function for (the closed $L_{\Phi}$-eigenvalue problem
if $\OM$ is a closed manifold) the Dirichlet $L_{\Phi}$-eigenvalue
problem. It is clear that $f$ is  an admissible function for the
respective $\triangle$-eigenvalue problem. Writing $\grad f (x) =
  \sum_{i=1}^{n}e_{i}(f)e_{i}(x)$ we have that \begin{eqnarray}\label{Obs3}\vert \Phi^{1/2}\,\grad f \vert^{2}(x)  &
  =&
   \langle \Phi \,\grad f, \,\grad f\rangle (x) \nonumber \\
   & = &  \langle \sum_{i=1}^{n}\mu_{i}(x)e_{i}(f)e_{i},
    \,\sum_{i=1}^{n}e_{i}(f)e_{i}\rangle\\
     & = & \sum_{i=1}^{n}\mu_{i}(x)e_{i}(f)^{2}(x).\nonumber
     \end{eqnarray}From (\ref{Obs3}) we have that
     \begin{equation}\label{Obs4}\nu (\Phi, M)\cdot \vert \grad f\vert^{2}(x)\geq\vert \Phi^{1/2}\,\grad f
    \vert^{2}(x)\geq \mu (\Phi, M)\cdot
  \vert \grad f\vert^{2}(x)
     \end{equation}and
     \begin{equation}\label{Obs5}\nu (\Phi, M)\cdot \displaystyle\frac{\int_{M}\vert\grad f\vert^{2}}{\int_{M}f^{2}} \geq
     \displaystyle\frac{\int_{M}\vert \Phi^{1/2}\,\grad f
    \vert^{2}}{\int_{M}f^{2}}\geq \mu (\Phi, M)\cdot \displaystyle\frac{\int_{M}\vert \grad f\vert^{2}}{\int_{M}f^{2}}
     \end{equation}Taking the infimum over all admissible functions in (\ref{Obs5}) we
     obtain (\ref{Obs-eq1}).

 \end{document}